\newcommand{\xBc}{\langle}
\newcommand{\xBe}{\rangle}

\newcommand{\xbS}{\Sigma}

\newcommand{\xba}{\alpha}
\newcommand{\xbb}{\beta}

\newcommand{\xbd}{\delta}
\newcommand{\xbe}{\in}
\newcommand{\xbf}{\phi}

\newcommand{\xbm}{\mu}

\newcommand{\xbq}{\psi}

\newcommand{\xbs}{\sigma}

\newcommand{\xCK}{\times}

\newcommand{\xCN}{\neg}
\newcommand{\xCQ}{\emptyset}

\newcommand{\xCf}{\hspace{0.1em}}

\newcommand{\xCq}{\sim}

\newcommand{\xcA}{\forall}

\newcommand{\xcE}{\exists}

\newcommand{\xcH}{\not\Rightarrow}

\newcommand{\xcM}{\not\models}

\newcommand{\xcS}{\bigcap}

\newcommand{\xcU}{\bigwedge}
\newcommand{\xcV}{\bigcup}

\newcommand{\xcb}{\subset}
\newcommand{\xcc}{\subseteq}

\newcommand{\xce}{\not\in}

\newcommand{\xcg}{\geq}
\newcommand{\xch}{\Rightarrow}

\newcommand{\xcj}{\Leftrightarrow}
\newcommand{\xck}{\leq}

\newcommand{\xcm}{\models}
\newcommand{\xcn}{\hspace{0.2em}\sim\hspace{-0.9em}\mid\hspace{0.58em}}

\newcommand{\xco}{\vee}
\newcommand{\xcp}{\rightarrow}

\newcommand{\xcs}{\cap}
\newcommand{\xcu}{\wedge}
\newcommand{\xcv}{\cup}

\newcommand{\xcz}{\Box}

\newcommand{\xDC}{\hspace{2em}}
\newcommand{\xDH}{\item }

\newcommand{\xDc}{\ll}

\newcommand{\xdC}{\mbox{\boldmath$C$}}

\newcommand{\xdd}{{\cal D}}

\newcommand{\xdf}{{\cal F}}

\newcommand{\xdl}{{\cal L}}

\newcommand{\xdo}{{\cal O}}
\newcommand{\xdp}{{\cal P}}

\newcommand{\xEI}{\begin{itemize}}
\newcommand{\xEJ}{\end{itemize}}

\newcommand{\xEd}{\neq}
\newcommand{\xEh}{\begin{enumerate}}
\newcommand{\xEj}{\end{enumerate}}

\newcommand{\xeA}{\nabla}
\newcommand{\xeB}{\not\prec}
\newcommand{\xeC}{\not\preceq}

\newcommand{\xeY}{\triangle}

\newcommand{\xeb}{\prec}
\newcommand{\xec}{\preceq}

\newcommand{\xFO}{\parallel}

\newcommand{\bl}{\begin{lemma} \rm}
\newcommand{\el}{\end{lemma}}
\newcommand{\br}{\begin{remark} \rm}
\newcommand{\er}{\end{remark}}
\newcommand{\be}{\begin{example} \rm}
\newcommand{\ee}{\end{example}}
\newcommand{\bco}{\begin{corollary} \rm}
\newcommand{\eco}{\end{corollary}}
\newcommand{\bc}{\begin{claim} \rm}
\newcommand{\ec}{\end{claim}}
\newcommand{\bfa}{\begin{fact} \rm}
\newcommand{\efa}{\end{fact}}
\newcommand{\bp}{\begin{proposition} \rm}
\newcommand{\ep}{\end{proposition}}
\newcommand{\bd}{\begin{definition} \rm}
\newcommand{\ed}{\end{definition}}
\newcommand{\bcs}{\begin{construction} \rm}
\newcommand{\ecs}{\end{construction}}
\newcommand{\bcd}{\begin{condition} \rm}
\newcommand{\ecd}{\end{condition}}
\newcommand{\bt}{\begin{theorem} \rm}
\newcommand{\et}{\end{theorem}}
\newcommand{\bn}{\begin{notation} \rm}
\newcommand{\en}{\end{notation}}
\newcommand{\bfi}{\begin{bild} \rm}
\newcommand{\efi}{\end{bild}}
\newcommand{\bsta}{\begin{statement} \rm}
\newcommand{\esta}{\end{statement}}
\newcommand{\bcom}{\begin{comment} \rm}
\newcommand{\ecom}{\end{comment}}
\newcommand{\bdia}{\begin{diagram} \rm}
\newcommand{\edia}{\end{diagram}}

\newcommand{\bfc}{\begin{figure}[htb] \begin{center}}
\newcommand{\efc}{\end{center} \end{figure}}

\documentclass{article}

\usepackage{amssymb,latexsym,epic,eepic,rotating}

\title{
A semantics for obligations -
\\
Local and global properties of obligations
\\
\thanks{
Paper 332
}
}

\author{Dov M Gabbay
\thanks{
Dov.Gabbay@kcl.ac.uk, www.dcs.kcl.ac.uk/staff/dg
} \\
King's College, London
\thanks{
Department of Computer Science, King's College London, Strand,
London WC2R 2LS, UK
} \\ \\
Karl Schlechta
\thanks{
ks@cmi.univ-mrs.fr, karl.schlechta@web.de, http://www.cmi.univ-mrs.fr/ $\sim$ ks
} \\
Laboratoire d'Informatique Fondamentale de Marseille
\thanks{
UMR 6166, CNRS and Universit\'{e} de Provence,
Address: CMI, 39, rue Joliot-Curie, F-13453 Marseille Cedex 13, France
}
}

\begin{document}

\newtheorem{lemma}{Lemma}[section]
\newtheorem{theorem}[lemma]{Theorem}
\newtheorem{proposition}[lemma]{Proposition}
\newtheorem{corollary}[lemma]{Corollary}
\newtheorem{claim}[lemma]{Claim}
\newtheorem{fact}[lemma]{Fact}
\newtheorem{remark}[lemma]{Remark}
\newtheorem{definition}{Definition}[section]
\newtheorem{construction}{Construction}[section]
\newtheorem{condition}{Condition}[section]
\newtheorem{example}{Example}[section]
\newtheorem{notation}{Notation}[section]
\newtheorem{bild}{Figure}[section]
\newtheorem{comment}{Comment}[section]
\newtheorem{statement}{Statement}[section]
\newtheorem{diagram}{Diagram}[section]

\renewcommand{\labelenumi}
  {(\arabic{enumi})}
\renewcommand{\labelenumii}
  {(\arabic{enumi}.\arabic{enumii})}
\renewcommand{\labelenumiii}
  {(\arabic{enumi}.\arabic{enumii}.\arabic{enumiii})}
\renewcommand{\labelenumiv}
  {(\arabic{enumi}.\arabic{enumii}.\arabic{enumiii}.\arabic{enumiv})}

\maketitle

\begin{abstract}

We analyze a number of properties obligations have or should have.
We will not present a definite answer what an obligation ``is", but
rather point out the numerous possibilities to consider, as well as some
of their interrelations.

\end{abstract}

\tableofcontents

\setcounter{secnumdepth}{3}
\setcounter{tocdepth}{3}

%
%
%
\section{
Introductory remarks
}

We see some relation of ``better'' as central for obligations.
Obligations determine what is ``better'' and what
is ``worse'', conversely, given such a relation of
``better'', we can define obligations.
The problems lie, in our opinion, in the fact that an adequate treatment
of such a relation is somewhat complicated, and leads to many
ramifications.

On the other hand, obligations have sufficiently many things
in common with facts so we can in a useful way say that an
obligation is satisfied in a situation, and one can also define a notion
of derivation for obligations.

Our approach is almost exclusively semantical.
\subsection{
Context
}

The problem with formalisation using logic is that the natural
movements in the application area being formalised do not exactly
correspond to natural movements in the logic being used as a tool of
formalisation. Put differently, we may be able to express statement A
of the application area by a formula $ \xbf $ of the logic, but the
subtleties of the way A is manipulated in the application area cannot
be matched in the formal logic used.
This gives rise to paradoxes. To resolve the paradoxes one needs to
improve the logic.
So the progress in the formalisation program depends on the state of
development of logic itself.
Recent serious advances in logical tools by the authors of this paper
enable us to offer some better formalisations possibilities for the
notion of obligation. This is what we offer in this paper.

Historically, articles on Deontic Logic include collections of
problems, see e.g.  \cite{MDW94}, semantical
approaches, see e.g.  \cite{Han69},
and others, like  \cite{CJ02}.

Our basic idea is to see obligations as tightly connected to some
relation of ``better''. An immediate consequence is that negation, which
inverses such a relation, behaves differently in the case of obligations
and of classical logic. $('' And'' $ and ``or'' seem to show analogue
behaviour in
both logics.)
The relation of ``better'' has to be treated with some caution, however,
and we introduce and investigate local and global properties about
``better'' of
obligations. Most of these properties coincide in sufficiently nice
situations, in others, they are different.

We do not come to a final conclusion which properties obligations should
or should not have, perhaps this will be answered in future, perhaps there
is no
universal answer. We provide a list of ideas which seem reasonable to us.

Throughout, we work in a finite (propositional) setting.
\subsection{
Central idea
}

We see a tight connection between obligations and a relation of
``morally'' better between situations. Obligations are there to
guide us for ``better'' actions, and, conversely, given some
relation of ``better'', we can define obligations.

The problems lie in the fact that a simple approach via quality is not
satisfactory. We examine a number of somewhat more subtle ideas, some
of them also using a notion of distance.
\subsection{
A common property of facts and obligations
}

We are fully aware that an obligation has a conceptually different status
than a fact. The latter is true or false, an obligation has been set by
some
instance as a standard for behaviour, or whatever.

Still, we will say that an obligation holds in a situation, or that a
situation
satisfies an obligation. If the letter is in the mail box, the obligation
to
post it is satisfied, if it is in the trash bin, the obligation is not
satisfied. In some set of worlds, the obligation is satisfied, in the
complement
of this set, it is not satisfied. Thus, obligations behave in this respect
like facts, and we put for this aspect all distinctions between facts and
obligations aside.

Thus, we will treat obligations most of the time as subsets of the set of
all
models, but also sometimes as formulas. As we work mostly in a finite
propositional setting, both are interchangeable.

We are $ \xCf not$ concerned here with actions to fulfill obligations,
developments
or so, just simply situations which satisfy or not obligations.

This is perhaps also the right place to point out that one has to
distinguish
between facts that $ \xCf hold$ in ``good'' situations (they will be closed
under
arbitrary right weakening), and obligations which describe what $ \xCf
should$
be, they will not be closed under arbitrary right weakening. This article
is only about the latter.
\subsection{
Derivations of obligations
}

Again, we are aware that ``deriving'' obligations is different from
``deriving''
facts. Derivation of facts is supposed to conclude from truth to truth,
deriving obligations will be about concluding what can also be considered
an obligation, given some set of ``basic'' obligations. The parallel is
sufficiently strong to justify the double use of the word ``derive''.

Very roughly, we will say that conjunctions (or intersections) and
disjunctions (unions) of obligations lead in a reasonable way to derived
obligations, but negations do not. We take the Ross paradox (see below)
very seriously, it was, as a matter of fact, our starting point to avoid
it in a reasonable notion of derivation.

We mention two simple postulates derived obligations should probably
satisfy.

 \xEh

 \xDH Every original obligation should also be a derived obligation,
corresponding to $ \xba, \xbb \xcn \xba.$

 \xDH A derived obligation should not be trivial, i.e. neither empty nor
$U,$
the universe we work in.

 \xEj

The last property is not very important from an algebraic point of
view, and easily satisfiable, so we will not give it too much attention.
\subsection{
Orderings and obligations
}

There is, in our opinion, a deep connection between obligations and
orderings (and, in addition, distances), which works both ways.

First, given a set of obligations, we can say that one situation is
``better'' than a seond situation, if the first satisfies ``more'' obligations
than the second does. ``More'' can be measured by the set of obligations
satisfied, and then by the subset/superset relation, or by the number of
obligations.
Both are variants of the Hamming distance. ``Distance'' between two
situations can be measured by the set or number of obligations in which
they differ (i.e. one situation satisfies them, the other not).
In both cases, we will call the variants the set or the counting variant.

This is also the deeper reason why we have to be careful with negation.
Negation inverses such orderings, if $ \xbf $ is better than $ \xCN \xbf
,$ then $ \xCN \xbf $ is
worse than $ \xCN \xCN \xbf = \xbf.$ But in some reasonable sense $ \xcu
$ and $ \xco $ preserve the
ordering, thus they are compatibel with obligations.

Conversely, given a relation (of quality), we might for instance require
that obligations are closed under improvement. More subtle requirements
might work with distances.
The relations of quality and distance can be given abstractly (as the
notion of size used for ``soft'' obligations), or as above by a starting set
of
obligations. We will also define important
auxiliary concepts on such abstract relations.
\subsection{
Derivation revisited
}

A set of ``basic'' obligations generates an ordering and a distance between
situations, ordering and distance can be used to define properties
obligations should have. It is thus natural to define obligations derived
from the basic set as those sets of situations which satisfy the
desirable properties of obligations defined via the order and distance
generated
by the basic set of obligations.
Our derivation is thus a two step procedure: first generate the order and
distance, which define suitable sets of situations.

We will call properties which are defined without using distances global
properties (like closure under improving quality), properties involving
distance (like being a neighbourhood) local properties.
\subsection{
Relativization
}

An important problem is relativization. Suppose $ \xdo $ is a set of
obligations for
all possible situations, e.g. $O$ is the obligation to post the letter,
and
$O' $ is the obligation to water the flowers. Ideally, we will do both.
Suppose
we consider now a subset, where we cannot do both (e.g. for lack of time).
What are our obligations in this subset? Are they just the restrictions to
the subset? Conversely, if $O$ is an obligation for a subset of all
situations, is then some $O' $ with $O \xcc O' $ an obligation for the set
of all
situations?

In more complicated requirements, it might be reasonable e.g. to choose
the ideal situations still in the big set, even if they are not in the
subset to be considered, but use an approximation inside the subset.
Thus, relativizations present a non-trivial problem with many possible
solutions, and it seems doubtful whether a universally adequate solution
can be
found.
\subsection{
Numerous possibilities
}

Seeing the possibilities presented so far (set or counting order,
set or counting distance, various relativizations), we can already guess
that there are numerous possible reasonable approaches to what an
obligation
is or should be. Consequently, it seems quite impossible to pursue all
these combinations in detail. Thus, we concentrate mostly on one
combination, and leave it to the reader to fill in details for the others,
if (s)he is so interested.

We will also treat the defeasible case here. Perhaps somewhat
surprisingly,
this is straightforward, and largely due to the fact that there is one
natural
definition of ``big'' sets for the product set, given that ``big'' sets are
defined for the components. So there are no new possibilities to deal with
here.

The multitude of possibilities is somewhat disappointing. It may,
of course, be due to an incapacity of the present authors to find $ \xCf
the$
right notion. But it may also be a genuine problem of ambiguous
intuitions,
and thus generate conflicting views and misunderstandings on the one side,
and loopholes for not quite honest argumentation in practical juridical
reasoning on the other hand.
\subsection{
Notation
}

$ \xdp (X)$ will be the powerset of $X,$ $A \xcc B$ will mean that A is a
subset of $B,$ or
$A=B,$ $A \xcb B$ that A is a proper subset of $B.$
\subsection{
Overview of the article
}

We will work in a finite propositional setting, so there is a trivial
and 1-1 correspondence between formulas and model sets. Thus, we can
just work with model sets - which implies, of course, that obligations
will be robust under logical reformulation. So we will formulate most
results only for sets.

 \xEI

 \xDH In Section \ref{Section Definitions} (page \pageref{Section Definitions})
,
we give the basic definitions, together with
some simple results about those definitions.

 \xDH Section \ref{Section Phil} (page \pageref{Section Phil})  will
present a more philosophical discussion, with more
examples, and we will argue that our definitions are relevant for our
purpose.

 \xDH As said already, there seems to be a multitude of possible and
reasonable definitions of what an obligation can or should be, so we
limit our formal investigation to a few cases, this is given in
Section \ref{Section Exam} (page \pageref{Section Exam}).

 \xDH In Section \ref{Section What} (page \pageref{Section What}), we give
a tentative definition of an obligation.

 \xEJ

(The concept of
neighbourhood semantics is not new, and was already
introduced by D.Scott,  \cite{Sco70}, and R.Montague,  \cite{Mon70}. Further
investigations showed that it was also used by O.Pacheco,  \cite{Pac07},
precisely to avoid unwanted weakening for obligations. We came the other
way
and started with the concept of independent strengthening, see
Definition \ref{Definition Delta-O} (page \pageref{Definition Delta-O}), and
introduced the abstract
concept of neighbourhood semantics only at the end. This is one of the
reasons we also have different descriptions which turned out to be
equivalent:
we came from elsewhere.)
\section{
Basic definitions
}
\label{Section Definitions}

We give here all definitions needed for our treatment of obligations.
The reader may continue immediately
to Section \ref{Section Phil} (page \pageref{Section Phil}), and come
back to the present section whenever necessary.

Intuitively, $U$ is supposed to be the set of
models of some propositional language $ \xdl,$ but we will stay purely
algebraic
whenever possible. $U' \xcc U$ is some subset.

For ease of writing, we will here and later sometimes work with
propositional
variables, and also identify models with the formula describing them,
still
this is just shorthand for arbitrary sets and elements. pq will stand for
$p \xcu q,$ etc.

If a set $ \xdo $ of obligations is given, these will be just arbitrary
subsets
of the universe $U.$ We will also say that $ \xdo $ is over $U.$

Before we deepen the discussion of more conceptual aspects, we give some
basic
definitions (for which we claim no originality). We will need them quite
early,
and think it is better to put them together here, not to be read
immediately,
but so the reader can leaf back and find them easily.

We work here with a notion of size (for the defeasible case), a notion $d$
of
distance, and a quality relation $ \xec.$ The latter two can be given
abstractly, but
may also be defined from a set of (basic) obligations.

We use these notions to describe properties obligations should, in our
opinion,
have. A careful analysis will show later interdependencies between the
different
properties.
\subsection{
Size
}

For each $U' \xcc U$ we suppose an abstract notion of size to be given.
We may assume this notion to be a filter or an ideal. Coherence properties
between the filters/ideals of different $U',$ $U'' $ will be left open,
the reader may assume them to be the conditions of the system $P$ of
preferential
logic, see  \cite{GS09a}.

Given such notions of size on $U' $ and $U'',$ we will also need a notion
of size
on $U' \xCK U''.$ We take the evident solution:

\bd

$\hspace{0.01em}$


\label{Definition Size-Product}

Let a notion of ``big subset'' be defined by a principal filter for all $X
\xcc U$ and
all $X' \xcc U'.$ Thus, for all $X \xcc U$ there exists a fixed principal
filter
$ \xdf (X) \xcc \xdp (X),$ and likewise for all $X' \xcc U'.$ (This is
the situation in the case
of preferential structures, where $ \xdf (X)$ is generated by $ \xbm (X),$
the set of
minimal elements of $X.)$

Define now $ \xdf (X \xCK X' )$ as generated by $\{A \xCK A':$ $A \xbe
\xdf (X),$ $A' \xbe \xdf (X' )\},$ i.e.
if A is the smallest element of $ \xdf (X),$ $A' $ the smallest element of
$ \xdf (X' ),$
then $ \xdf (X \xCK X' ):=\{B \xcc X \xCK X':$ $A \xCK A' \xcc B\}.$

\ed

\bfa

$\hspace{0.01em}$


\label{Fact Product-Pref}

If $ \xdf (X)$ and $ \xdf (X' )$ are generated by preferential structures
$ \xec_{X},$ $ \xec_{X' },$
then $ \xdf (X \xCK X' )$ is
generated by the product structure defined by

$ \xBc x,x'  \xBe  \xec_{X \xCK X' } \xBc y,y'  \xBe $ $: \xcj $ $x \xec_{X}y$
and $x'
\xec_{X' }y'.$

\efa

\subparagraph{
Proof
}

$\hspace{0.01em}$


We will omit the indices of the orderings when this causes no confusion.

Let $A \xbe \xdf (X),$ $A' \xbe \xdf (X' ),$ i.e. A minimizes $X,$ $A' $
minimizes $X'.$ Let
$ \xBc x,x'  \xBe  \xbe X \xCK X',$ then there are $a \xbe A,$ $a' \xbe A' $
with
$a \xec x,$ $a' \xec x',$ so
$ \xBc a,a'  \xBe  \xec  \xBc x,x'  \xBe.$

Conversely, suppose $U \xcc X \xCK X',$ $U$ minimizes $X \xCK X'.$ but
there is no $A \xCK A' \xcc U$
s.t. $A \xbe \xdf (X),$ $A' \xbe \xdf (X' ).$ Assume $A= \xbm (X),$ $A' =
\xbm (X' ),$ so there is
$ \xBc a,a'  \xBe  \xbe \xbm (X) \xCK \xbm (X' ),$ $ \xBc a,a'  \xBe  \xce U.$
But only
$ \xBc a,a'  \xBe  \xec  \xBc a,a'  \xBe,$ and $U$ does not
minimize $X \xCK X',$ $contradiction.$

$ \xcz $
\\[3ex]

Note that a natural modification of our definition:

There is $A \xbe \xdf (X)$ s.t.
for all $a \xbe A$ there is a (maybe varying) $A'_{a} \xbe \xdf (X' ),$
and
$U:=\{ \xBc a,a'  \xBe :$ $a \xbe A,$ $a' \xbe A'_{a}\}$ as generating sets

will result in the same
definition, as our filters are principal, and thus stable under arbitrary
intersections.
\subsection{
Distance
}

We consider a set of sequences $ \xbS,$ for $x \xbe \xbS $ $x:I \xcp S,$
$I$ a finite index set, $S$
some set.
Often, $S$ will be $\{0,1\},$ $x(i)=1$ will mean that $x \xbe i,$ when $I
\xcc \xdp (U)$ and $x \xbe U.$
For abbreviation, we will call this (unsystematically, often context will
tell) the $ \xbe -$case.
Often, $I$ will be written $ \xdo,$ intuitively, $O \xbe \xdo $ is then
an obligation, and
$x(O)=1$ means $x \xbe O,$ or $x$ ``satisfies'' the obligation $O.$

\bd

$\hspace{0.01em}$


\label{Definition O-x}

In the $ \xbe -$case, set $ \xdo (x):=\{O \xbe \xdo:x \xbe O\}.$

\ed

\bd

$\hspace{0.01em}$


\label{Definition Distance}

Given $x,y \xbe \xbS,$ the Hamming distance comes in two flavours:

$d_{s}(x,y):=\{i \xbe I:x(i) \xEd y(i)\},$ the set variant,

$d_{c}(x,y):=card(d_{s}(x,y)),$ the counting variant.

We define $d_{s}(x,y) \xck d_{s}(x',y' )$ iff $d_{s}(x,y) \xcc d_{s}(x'
,y' ),$

thus, $s-$distances are not always comparabel.

\ed

\bfa

$\hspace{0.01em}$


\label{Fact Distance}

(1) In the $ \xbe -$case, we have $d_{s}(x,y)= \xdo (x) \xeY \xdo (y),$
where $ \xeY $ is the symmetric
set difference.

(2) $d_{c}$ has the normal addition, set union takes the role of addition
for $d_{s},$
$ \xCQ $ takes the role of 0 for $d_{s},$
both are distances in the following sense:

(2.1) $d(x,y)=0$ if $x=y,$ but not conversely,

(2.2) $d(x,y)=d(y,x),$

(2.3) the triangle inequality holds, for the set variant in the form
$d_{s}(x,z) \xcc d_{s}(x,y) \xcv d_{s}(y,z).$

(If $d(x,y)=0$ $ \xcH $ $x=y$ poses a problem, one can always consider
equivalence
classes.)

\efa

\subparagraph{
Proof
}

$\hspace{0.01em}$


(2.1) Suppose $U=\{x,y\},$ $ \xdo =\{U\},$ then $ \xdo (X)= \xdo (Y),$ but
$x \xEd y.$

(2.3) If $i \xce d_{s}(x,y) \xcv d_{s}(y,z),$ then $x(i)=y(i)=z(i),$ so
$x(i)=z(i)$ and
$i \xce d_{s}(x,z).$

The others are trivial.

$ \xcz $
\\[3ex]

\bd

$\hspace{0.01em}$


\label{Definition Between}

(1) We can define for any distance $d$ with some minimal requirements a
notion of
``between''.

If the codomain of $d$ has an ordering $ \xck,$ but no addition, we
define:

$ \xBc x,y,z \xBe _{d}$ $: \xcj $ $d(x,y) \xck d(x,z)$ and $d(y,z) \xck d(x,z).$

If the codomain has a commutative addition, we define

$ \xBc x,y,z \xBe _{d}$ $: \xcj $ $d(x,z)=d(x,y)+d(y,z)$ - in $d_{s}$ $+$ will
be
replaced by $ \xcv,$ i.e.

$ \xBc x,y,z \xBe _{s}$ $: \xcj $ $d(x,z)=d(x,y) \xcv d(y,z).$

For above two Hamming distances, we will write $ \xBc x,y,z \xBe _{s}$ and
$ \xBc x,y,z \xBe _{c}.$

(2) We further define:

$[x,z]_{d}:=\{y \xbe X: \xBc x,y,x \xBe _{d}\}$ - where $X$ is the set we work
in.

We will write $[x,z]_{s}$ and $[x,z]_{c}$ when appropriate.

(3) For $x \xbe U,$ $X \xcc U$ set $x \xFO_{d}X$ $:=$ $\{x' \xbe X: \xCN
\xcE x'' \xEd x' \xbe X.d(x,x' ) \xcg d(x,x'' )\}.$

Note that, if $X \xEd \xCQ,$ then $x \xFO X \xEd \xCQ.$

We omit the index when this does not cause confusion. Again, when
adequate,
we write $ \xFO_{s}$ and $ \xFO_{c}.$

\ed

For problems with characterizing ``between'' see
 \cite{Sch04}.

\bfa

$\hspace{0.01em}$


\label{Fact Between}

(0) $ \xBc x,y,z \xBe _{d}$ $ \xcj $ $ \xBc z,y,x \xBe _{d}.$

Consider the situation of a set of sequences $ \xbS.$

Let $A:=A_{ \xbs, \xbs '' }:=\{ \xbs ': \xcA i \xbe I( \xbs (i)= \xbs ''
(i) \xcp \xbs ' (i)= \xbs (i)= \xbs '' (i))\}.$
Then

(1) If $ \xbs ' \xbe A,$ then $d_{s}( \xbs, \xbs '' )=d_{s}( \xbs, \xbs
' ) \xcv d_{s}( \xbs ', \xbs '' ),$ so $ \xBc  \xbs, \xbs ', \xbs ''  \xBe
_{s}.$

(2) If $ \xbs ' \xbe A$ and $S$ consists of 2 elements (as in classical
2-valued
logic), then $d_{s}( \xbs, \xbs ' )$ and $d_{s}( \xbs ', \xbs '' )$ are
disjoint.

(3) $[ \xbs, \xbs '' ]_{s}=A.$

(4) If, in addition, $S$ consists of 2 elements, then $[ \xbs, \xbs ''
]_{c}=A.$

\efa

\subparagraph{
Proof
}

$\hspace{0.01em}$


(0) Trivial.

(1) `` $ \xcc $ '' follows from Fact \ref{Fact Distance} (page \pageref{Fact
Distance}), (2.3).

Conversely, if e.g. $i \xbe d_{s}( \xbs, \xbs ' ),$ then
by prerequisite $i \xbe d_{s}( \xbs, \xbs '' ).$

(2) Let $i \xbe d_{s}( \xbs, \xbs ' ) \xcs d_{s}( \xbs ', \xbs '' ),$
then $ \xbs (i) \xEd \xbs ' (i)$ and $ \xbs ' (i) \xEd \xbs '' (i),$
but then by $card(S)=2$ $ \xbs (i)= \xbs '' (i),$ but $ \xbs ' \xbe A,$
$contradiction.$

We turn to (3) and (4):

If $ \xbs ' \xce A,$ then there is $i' $ s.t. $ \xbs (i' )= \xbs '' (i' )
\xEd \xbs ' (i' ).$ On the other hand,
for all $i$ s.t. $ \xbs (i) \xEd \xbs '' (i)$ $i \xbe d_{s}( \xbs, \xbs '
) \xcv d_{s}( \xbs ', \xbs '' ).$ Thus:

(3) By (1) $ \xbs ' \xbe A$ $ \xch $ $ \xBc  \xbs, \xbs ', \xbs ''  \xBe _{s}.$
Suppose $ \xbs ' \xce A,$ so there is $i' $ s.t.
$i' \xbe d_{s}( \xbs, \xbs ' )-d_{s}( \xbs, \xbs '' ),$ so $< \xbs,
\xbs ', \xbs '' >_{s}$ cannot be.

(4) By (1) and (2) $ \xbs ' \xbe A$ $ \xch $ $< \xbs, \xbs ', \xbs ''
>_{c}.$ Conversely, if $ \xbs ' \xce A,$ then
$card(d_{s}( \xbs, \xbs ' ))+card(d_{s}( \xbs ', \xbs '' )) \xcg
card(d_{s}( \xbs, \xbs '' ))+2.$

$ \xcz $
\\[3ex]

\bd

$\hspace{0.01em}$


\label{Definition L-And}

Given a finite propositional laguage $ \xdl $ defined by the set $v( \xdl
)$ of
propositional
variables, let $ \xdl_{ \xcu }$ be the set of all consistent conjunctions
of
elements from $v( \xdl )$ or their negations. Thus, $p \xcu \xCN q \xbe
\xdl_{ \xcu }$ if $p,q \xbe v( \xdl ),$ but
$p \xco q,$ $ \xCN (p \xcu q) \xce \xdl_{ \xcu }.$ Finally, let $ \xdl_{
\xco \xcu }$ be the set of all (finite)
disjunctions of formulas from $ \xdl_{ \xcu }.$ (As we will later not
consider all
formulas from $ \xdl_{ \xcu },$ this will be a real restriction.)

Given a set of models $M$ for a finite language $ \xdl,$ define
$ \xbf_{M}:= \xcU \{p \xbe v( \xdl ): \xcA m \xbe M.m(p)=v\} \xcu \xcU \{
\xCN p:p \xbe v( \xdl ), \xcA m \xbe M.m(p)=f\} \xbe \xdl_{ \xcu }.$
(If there are no such $p,$ set $ \xbf_{M}:=TRUE.)$

This is the strongest $ \xbf \xbe \xdl_{ \xcu }$ which holds in $M.$

\ed

\bfa

$\hspace{0.01em}$


\label{Fact Hamming-Neighbourhood}

If $x,y$ are models, then $[x,y]=M( \xbf_{\{x,y\}}).$ $ \xcz $
\\[3ex]

\efa

\subparagraph{
Proof
}

$\hspace{0.01em}$


$m \xbe [x,y]$ $ \xcj $ $ \xcA p(x \xcm p,y \xcm p \xch m \xcm p$ and $x
\xcM p,y \xcM p \xch m \xcM p),$
$m \xcm \xbf_{\{x,y\}}$ $ \xcj $ $m \xcm \xcU \{p:x(p)=y(p)=v\} \xcu \xcU
\{ \xCN p:x(p)=y(p)=f\}.$
\subsection{
Quality and closure
}

\bd

$\hspace{0.01em}$


\label{Definition Closed}

Given any relation $ \xec $ (of quality), we say that $X \xcc U$ is
(downward) closed (with
respect to $ \xec )$ iff $ \xcA x \xbe X \xcA y \xbe U(y \xec x$ $ \xch $
$y \xbe X).$

\ed

(Warning, we follow the preferential tradition, ``smaller'' will
mean ``better''.)

We think that being closed is a desirable property for obligations: what
is
at least as good as one element in the obligation should be ``in'', too.

\bfa

$\hspace{0.01em}$


\label{Fact Subset-Closure}

Let $ \xec $ be given.

(1) Let $D \xcc U' \xcc U'',$ $D$ closed in $U'',$ then $D$ is also
closed in $U'.$

(2) Let $D \xcc U' \xcc U'',$ $D$ closed in $U',$ $U' $ closed in $U''
,$ then $D$ is
closed in $U''.$

(3) Let $D_{i} \xcc U' $ be closed for all $i \xbe I,$ then so are $ \xcV
\{D_{i}:i \xbe I\}$ and $ \xcS \{D_{i}:i \xbe I\}.$

\efa

\subparagraph{
Proof
}

$\hspace{0.01em}$


(1) Trivial.

(2) Let $x \xbe D \xcc U',$ $x' \xec x,$ $x' \xbe U'',$ then $x' \xbe U'
$ by closure of $U'',$ so $x' \xbe D$ by
closure of $U'.$

(3) Trivial.

$ \xcz $
\\[3ex]

We may have an abstract relation $ \xec $ of quality on the domain, but we
may also
define it from the structure of the sequences, as we will do now.

\bd

$\hspace{0.01em}$


\label{Definition Quality}

Consider the case of sequences.

Given a relation $ \xec $ (of quality) on the codomain, we extend this to
sequences
in $ \xbS:$

$x \xCq y$ $: \xcj $ $ \xcA i \xbe I(x(i) \xCq y(i))$

$x \xec y$ $: \xcj $ $ \xcA i \xbe I(x(i) \xec y(i))$

$x \xeb y$ $: \xcj $ $ \xcA i \xbe I(x(i) \xec y(i))$ and $ \xcE i \xbe
I(x(i) \xeb y(i))$

In the $ \xbe -$case, we will consider $x \xbe i$ better than $x \xce i.$
As we have only two
values, true/false, it is easy to count the positive and negative cases
(in more complicated situations, we might be able to multiply), so we have
an analogue of the two Hamming distances, which we might call the Hamming
quality relations.

Let $ \xdo $ be given now.

(Recall that we follow the preferential tradition, ``smaller'' will
mean ``better''.)

$x \xCq_{s}y$ $: \xcj $ $ \xdo (x)= \xdo (y),$

$x \xec_{s}y$ $: \xcj $ $ \xdo (y) \xcc \xdo (x),$

$x \xeb_{s}y$ $: \xcj $ $ \xdo (y) \xcb \xdo (x),$

$x \xCq_{c}y$ $: \xcj $ $card( \xdo (x))=card( \xdo (y)),$

$x \xec_{c}y$ $: \xcj $ $card( \xdo (y)) \xck card( \xdo (x)),$

$x \xeb_{c}y$ $: \xcj $ $card( \xdo (y))<card( \xdo (x)).$

\ed

The requirement of closure causes a problem for
the counting approach: Given e.g. two obligations $O,$ $O',$ then any two
elements
in just one obligation have the same quality, so if one is in, the other
should
be, too. But this prevents now any of the original obligations to have the
desirable property of closure. In the counting case, we will obtain a
ranked
structure, where elements satisfy
0, 1, 2, etc. obligations, and we are unable to differentiate inside those
layers. Moreover, the set variant seems to be closer to logic, where we do
not
count the propositional variables which hold in a model, but consider them
individually. For these reasons, we will not pursue the counting approach
as
systematically as the set approach. One should, however, keep in mind that
the
counting variant gives a ranking relation of quality, as all qualities are
comparable, and the set variant does not. A ranking seems to be
appreciated
sometimes in the literature, though we are not really sure why.

Of particular interest is the combination of $d_{s}$ and $ \xec_{s}$
$(d_{c}$ and $ \xec_{c})$
respectively - where by $ \xec_{s}$ we also mean $ \xeb_{s}$ and $
\xCq_{s},$ etc. We turn to
this now.

\bfa

$\hspace{0.01em}$


\label{Fact Quality-Distance}

We work in the $ \xbe -$case.

(1) $x \xec_{s}y$ $ \xch $ $d_{s}(x,y)= \xdo (x)- \xdo (y)$

Let $a \xeb_{s}b \xeb_{s}c.$ Then

(2) $d_{s}(a,b)$ and $d_{s}(b,c)$ are not comparable,

(3) $d_{s}(a,c)=d_{s}(a,b) \xcv d_{s}(b,c),$ and thus $b \xbe [a,c]_{s}.$

This does not hold in the counting variant, as Example \ref{Example Count} (page
\pageref{Example Count})  shows.

(4) Let $x \xeb_{s}y$ and $x' \xeb_{s}y$ with $x,x' $ $
\xeb_{s}-$incomparabel. Then $d_{s}(x,y)$ and $d_{s}(x',y)$
are incomparable.

(This does not hold in the counting variant, as then all distances are
comparable.)

(5) If $x \xeb_{s}z,$ then for all $y \xbe [x,z]_{s}$ $x \xec_{s}y
\xec_{s}z.$

\efa

\subparagraph{
Proof
}

$\hspace{0.01em}$


(1) Trivial.

(2) We have $ \xdo (c) \xcb \xdo (b) \xcb \xdo (a),$ so the results
follows from (1).

(3) By definition of $d_{s}$ and (1).

(4) $x$ and $x' $ are $ \xec_{s}-incomparable,$ so there are $O \xbe \xdo
(x)- \xdo (x' ),$
$O' \xbe \xdo (x' )- \xdo (x).$

As $x,x' \xeb_{s}y,$ $O,O' \xce \xdo (y),$ so $O \xbe d_{s}(x,y)-d_{s}(x'
,y),$
$O' \xbe d_{s}(x',y)-d_{s}(x,y).$

(5) $x \xeb_{s}z$ $ \xch $ $ \xdo (z) \xcb \xdo (x),$ $d_{s}(x,z)= \xdo
(x)- \xdo (z).$ By prerequisite
$d_{s}(x,z)=d_{s}(x,y) \xcv d_{s}(y,z).$
Suppose $x \xeC_{s}y.$ Then there is $i \xbe \xdo (y)- \xdo (x) \xcc
d_{s}(x,y),$ so
$i \xce \xdo (x)- \xdo (z)=d_{s}(x,z),$ $contradiction.$

Suppose $y \xeC_{s}z.$ Then there is $i \xbe \xdo (z)- \xdo (y) \xcc
d_{s}(y,z),$ so
$i \xce \xdo (x)- \xdo (z)=d_{s}(x,z),$ $contradiction.$

$ \xcz $
\\[3ex]

\be

$\hspace{0.01em}$


\label{Example Count}

In this and similar examples, we will use the model notation. Some
propositional variables $p,$ $q,$ etc. are given, and models are described
by
$p \xCN qr,$ etc. Moreover, the propositional variables are the
obligations, so
in this example we have the obligations $p,$ $q,$ $r.$

Consider $x:= \xCN p \xCN qr,$ $y:=pq \xCN r,$ $z:= \xCN p \xCN q \xCN r.$
Then $y \xeb_{c}x \xeb_{c}z,$ $d_{c}(x,y)=3,$
$d_{c}(x,z)=1,$ $d_{c}(z,y)=2,$ so $x \xce [y,z]_{c}.$ $ \xcz $
\\[3ex]

\ee

\bd

$\hspace{0.01em}$


\label{Definition Quality-Extension}

Given a quality relation $ \xeb $ between elements, and a distance $d,$ we
extend the
quality relation to sets and define:

(1) $x \xeb Y$ $: \xcj $ $ \xcA y \xbe (x \xFO Y).x \xeb y.$ (The closest
elements - i.e. there are no
closer ones - of $Y,$ seen from $x,$ are less good than $x.)$

analogously $X \xeb y$ $: \xcj $ $ \xcA x \xbe (y \xFO X).x \xeb y$

(2) $X \xeb_{l}Y$ $: \xcj $ $ \xcA x \xbe X.x \xeb Y$ and $ \xcA y \xbe
Y.X \xeb y$
(X is locally better than $Y).$

When necessary, we will write $ \xeb_{l,s}$ or $ \xeb_{l,c}$ to
distinguish the
set from the counting variant.

For the next definition, we use the notion of size: $ \xeA \xbf $ iff for
almost all $ \xbf $
holds i.e. the set of exceptions is small.

(3) $X \xDc_{l}Y$ $: \xcj $ $ \xeA x \xbe X.x \xeb Y$ and $ \xeA y \xbe
Y.X \xeb y.$

We will likewise write $ \xDc_{l,s}$ etc.

This definition is supposed to capture quality difference under minimal
change,
the ``ceteris paribus'' idea: $X \xeb_{l} \xdC X$ should hold for an
obligation $X.$
Minimal change is coded by $ \xFO,$ and ``ceteris paribus'' by minimal
change.

\ed

\bfa

$\hspace{0.01em}$


\label{Fact General-Obligation}

If $X \xeb_{l} \xdC X,$ and $x \xbe U$ an optimal point (there is no
better one), then $x \xbe X.$

\efa

\subparagraph{
Proof
}

$\hspace{0.01em}$


If not, then take $x' \xbe X$ closest to $x,$ this must be better than
$x,$ contradiction.
$ \xcz $
\\[3ex]

\bfa

$\hspace{0.01em}$


\label{Fact Local-Global}

Take the set version.

If $X \xeb_{l,s} \xdC X,$ then $X$ is downward $ \xeb_{s}-closed.$

\efa

\subparagraph{
Proof
}

$\hspace{0.01em}$


Suppose $X \xeb_{l,s} \xdC X,$ but $X$ is not downward closed.

Case 1: There are $x \xbe X,$ $y \xce X,$ $y \xCq_{s}x.$ Then $y \xbe x
\xFO_{s} \xdC X,$ but $x \xeB y,$ $contradiction.$

Case 2:
There are $x \xbe X,$ $y \xce X,$
$y \xeb_{s}x.$ By $X \xeb_{l,s} \xdC X,$ the elements in $X$ closest to
$y$ must be better than $y.$
Thus, there is $x' \xeb_{s}y,$ $x' \xbe X,$ with minimal distance from
$y.$ But then
$x' \xeb_{s}y \xeb_{s}x,$ so $d_{s}(x',y)$ and $d_{s}(y,x)$ are
incomparable by
Fact \ref{Fact Quality-Distance} (page \pageref{Fact Quality-Distance}), so $x$
is among those
with minimal distance from $y,$ so $X \xeb_{l,s} \xdC X$ does not hold. $
\xcz $
\\[3ex]

\be

$\hspace{0.01em}$


\label{Example Dependent-2}

We work with the set variant.

This example shows that $ \xec_{s}-$closed does not imply $X \xeb_{l,s}
\xdC X,$ even
if $X$ contains the best elements.

Let $ \xdo:=\{p,q,r,s\},$ $U':=\{x:=p \xCN q \xCN r \xCN s,$ $y:= \xCN
pq \xCN r \xCN s,$ $x':=pqrs\},$ $X:=\{x,x' \}.$
$x' $ is the best element of $U',$ so $X$ contains the best elements,
and $X$ is downward closed in $U',$
as $x$ and $y$ are not comparable. $d_{s}(x,y)=\{p,q\},$ $d_{s}(x'
,y)=\{p,r,s\},$ so the
distances from $y$ are not comparable, so $x$ is among the closest
elements in $X,$
seen from $y,$ but $x \xeB_{s}y.$

The lack of comparability is essential here, as the following Fact shows.

$ \xcz $
\\[3ex]

\ee

We have, however, for the counting variant:

\bfa

$\hspace{0.01em}$


\label{Fact Count-Closed}

Consider the counting variant. Then

If $X$ is downward closed, then $X \xeb_{l,c} \xdC X.$

\efa

\subparagraph{
Proof
}

$\hspace{0.01em}$


Take any $x \xbe X,$ $y \xce X.$ We have $y \xec_{c}x$ or $x \xeb_{c}y,$
as any two elements are
$ \xec_{c}-$comparabel. $y \xec_{c}x$ contradicts closure, so $x
\xeb_{c}y,$ and $X \xeb_{l,c} \xdC X$ holds
trivially. $ \xcz $
\\[3ex]
\subsection{
Neighbourhood
}

\bd

$\hspace{0.01em}$


\label{Definition Neighbourhood}

Given a distance $d,$ we define:

(1) Let $X \xcc Y \xcc U',$ then $Y$ is a neighbourhood of $X$ in $U' $
iff

$ \xcA y \xbe Y \xcA x \xbe X(x$ is closest to $y$ among all $x' $ with
$x' \xbe X$ $ \xch $ $[x,y] \xcs U' \xcc Y).$

(Closest means that there are no closer ones.)

When we also have a quality relation $ \xeb,$ we define:

(2) Let $X \xcc Y \xcc U',$ then $Y$ is an improving neighbourhood of $X$
in $U' $ iff

$ \xcA y \xbe Y \xcA x((x$ is closest to $y$ among all $x' $ with $x' \xbe
X$ and $x' \xec y)$ $ \xch $ $[x,y] \xcs U' \xcc Y).$

When necessary, we will have to say for (3) and (4) which variant, i.e.
set or counting, we mean.

\ed

\bfa

$\hspace{0.01em}$


\label{Fact Neighbourhood}

(1) If $X \xcc X' \xcc \xbS,$ and $d(x,y)=0$ $ \xch $ $x=y,$ then $X$ and
$X' $ are Hamming neighbourhoods
of $X$ in $X'.$

(2) If $X \xcc Y_{j} \xcc X' \xcc \xbS $ for $j \xbe J,$ and all $Y_{j}$
are Hamming Neighbourhoods
of $X$ in $X',$ then so are $ \xcV \{Y_{j}:j \xbe J\}$ and $ \xcS
\{Y_{j}:j \xbe J\}.$

\efa

\subparagraph{
Proof
}

$\hspace{0.01em}$


(1) is trivial (we need here that $d(x,y)=0$ $ \xch $ $x=y).$

(2) Trivial.

$ \xcz $
\\[3ex]
\subsection{
Unions of intersections and other definitions
}

\bd

$\hspace{0.01em}$


\label{Definition ui}

Let $ \xdo $ over $U$ be given.

$X \xcc U' $ is $ \xCf (ui)$ (for union of intersections) iff there is a
family $ \xdo_{i} \xcc \xdo,$
$i \xbe I$ s.t. $X=( \xcV \{ \xcS \xdo_{i}:i \xbe I\}) \xcs U'.$

\ed

Unfortunately, this definition is not very useful for simple
relativization.

\bd

$\hspace{0.01em}$


\label{Definition Validity}

Let $ \xdo $ be over $U.$ Let $ \xdo ' \xcc \xdo.$ Define for $m \xbe U$
and $ \xbd: \xdo ' \xcp 2=\{0,1\}$

$m \xcm \xbd $ $: \xcj $ $ \xcA O \xbe \xdo ' (m \xbe O \xcj \xbd (O)=1)$

\ed

\bd

$\hspace{0.01em}$


\label{Definition Independence}

Let $ \xdo $ be over $U.$

$ \xdo $ is independent iff $ \xcA \xbd: \xdo \xcp 2. \xcE m \xbe U.m
\xcm \xbd.$

\ed

Obviously, independence does not inherit downward to subsets of $U.$

\bd

$\hspace{0.01em}$


\label{Definition Delta-O}

This definition is only intended for the set variant.

Let $ \xdo $ be over $U.$

$ \xdd ( \xdo ):=\{X \xcc U':$ $ \xcA \xdo ' \xcc \xdo $ $ \xcA \xbd:
\xdo ' \xcp 2$

$ \xDC $ $(( \xcE m,m' \xbe U,$ $m,m' \xcm \xbd,$ $m \xbe X,m' \xce X)$ $
\xch $ $( \xcE m'' \xbe X.m'' \xcm \xbd \xcu m'' \xeb_{s}m' ))\}$

This property expresses that we can satisfy obligations independently: If
we
respect $O,$ we can, in addition,
respect $O',$ and if we are hopeless kleptomaniacs, we may still not be a
murderer. If $X \xbe \xdd ( \xdo ),$ we can go from $U-X$
into $X$ by improving on all $O \xbe \xdo,$ which we have not fixed by $
\xbd,$ if $ \xbd $ is
not too rigid.
\section{
Philosophical discussion of obligations
}
\label{Section Phil}

\ed

We take now a closer look at obligations, in particular at the
ramifications
of the treatment of the relation ``better''. Some aspects of obligations
will also
need a notion of distance, we call them local properties of obligations.
\subsection{
A fundamental difference between facts and obligations: asymmetry and
negation
}

There is an important difference between facts and obligations. A
situation
which satisfies an obligation is in some sense ``good'', a situation which
does not, is in some sense ``bad''. This is not true of facts. Being
``round'' is a priori not better than ``having corners'' or vice versa. But
given
the obligation to post the letter, the letter in the mail box is ``good'',
the
letter in he trash bin is ``bad''. Consequently, negation has to play
different role for obligations and for facts.

This is a fundamental property, which can also be
found in orders, planning (we move towards the goal or not), reasoning
with
utility (is $ \xbf $ or $ \xCN \xbf $ more useful?), and probably others,
like perhaps the
Black Raven paradox.

We also think that the Ross paradox (see below) is a true paradox, and
should be
avoided. A closer look shows that this paradox involves arbitrary
weakening,
in particular by the ``negation'' of an obligation. This was a starting
point of our analysis.

``Good'' and ``bad'' cannot mean that any situation satisfying obligation $O$
is better than any situation not satisfying $O,$ as the following example
shows.

\be

$\hspace{0.01em}$


\label{Example 3-Obligations}

If we have three independent
and equally strong obligations, $O,$ $O',$ $O'',$ then a situation
satisfying
$O$ but neither $O' $ nor $O'' $ will not be better than one satisfying
$O' $ and $O'',$
but not $O.$

\ee

We have to introduce some kind of ``ceteris paribus''. All other
things being equal, a situation satisfying $O$ is better than a situation
not satisfying $O,$
see Section \ref{Section Ceteris} (page \pageref{Section Ceteris}).

\be

$\hspace{0.01em}$


\label{Example Ross-Paradox}

The original version of the Ross paradox reads: If we have the obligation
to post the letter, then we have the obligation to post or burn the
letter.
Implicit here is the background knowledge that burning the letter implies
not
to post it, and is even worse than not posting it.

We prefer a modified version, which works with two independent
obligations:
We have the obligation to post the letter, and we have the obligation to
water
the plants. We conclude by unrestricted weakening that we have the
obligation to
post the letter or $ \xCf not$ to water the plants. This is obvious
nonsense.

It is not the ``or'' itself which is the problem. For instance, in case of
an
accident, to call an ambulance or to help the victims by giving first aid
is a
perfectly reasonable obligation. It is the negation of the obligation to
water the plants which is the problem. More generally,
it must not be that the system of suitable sets is closed under
arbitrary supersets,
otherwise we have closure under arbitrary right weakening, and thus the
Ross
paradox. Notions like ``big subset'' or ``small exception sets'' from the
semantics of nonmonotonic logics are closed under
supersets, so they are not suitable.
\subsection{
``And'' and ``or'' for obligations
}

\ee

``Not'' behaves differently for facts and for obligations. If $O$ and $O' $
are
obligations, can $O \xcu O' $ be considered an obligation? We think, yes.
``Ceteris paribus'', satisfying $O$ and $O' $ together is better than not to
do
so. If is the obligation to post the letter, $O' $ to water the plants,
then
doing both is good, and better than doing none, or only one. Is $O \xco O'
$
an obligation? Again, we think, yes. Satisfying one (or even both, a
non-exclusive or) is better than doing nothing. We might not have enough
time
to do both, so we do our best, and water the plants or post the letter.
Thus, if $ \xba $ and $ \xbb $ are obligations, then so will be $ \xba
\xcu \xbb $ and $ \xba \xco \xbb,$ but not
anything involving $ \xCN \xba $ or $ \xCN \xbb.$ (In a non-trivial
manner, leaving aside
tautologies and contradictions which have to be considered separately.)
To summarize: ``and'' and ``or'' preserve the asymmetry, ``not'' does not,
therefore
we can combine obligations using ``and'' and ``or'', but not ``not''. Thus,
a reasonable notion of derivation of obligations will work with $ \xcu $
and $ \xco,$ but
not with $ \xCN.$

We should not close under inverse $ \xcu,$ i.e. if $ \xbf \xcu \xbf ' $
is an obligation, we
should not conclude that $ \xbf $ and $ \xbf ' $ separately are
obligations, as the following
example shows.

\be

$\hspace{0.01em}$


\label{Example Not-Inv-And}

Let $p$ stand for: post letter, $w:$ water plants, $s:$ strangle
grandmother.

Consider now $ \xbf \xcu \xbf ',$ where $ \xbf =p \xco ( \xCN p \xcu \xCN
w),$ $ \xbf ' =p \xco ( \xCN p \xcu w \xcu s).$
$ \xbf \xcu \xbf ' $ is equivalent to $p$ - though it is perhaps a bizarre
way to express
the obligation to post the letter. $ \xbf $ leaves us the possibility not
to water the
plants, and $ \xbf ' $ to strangle the grandmother, and neither seem good
obligations.
$ \xcz $
\\[3ex]

\ee

\br

$\hspace{0.01em}$


\label{Remark Not-Inv-And}

This is particularly important in the case of soft obligations, as we see
now,
when we try to apply the rules of preferential reasoning to obligations.

One of the rules of preferential reasoning is the $ \xCf (OR)$ rule:

$ \xbf \xcn \xbq,$ $ \xbf ' \xcn \xbq $ $ \xch $ $ \xbf \xco \xbf ' \xcn
\xbq.$

Suppose we have $ \xbf \xcn \xbq ' \xcu \xbq '',$ and $ \xbf ' \xcn \xbq
'.$ We might be tempted to split
$ \xbq ' \xcu \xbq '' $ - as $ \xbq ' $ is a ``legal'' obligation, and
argue: $ \xbf \xcn \xbq ' \xcu \xbq '',$ so
$ \xbf \xcn \xbq ',$ moreover $ \xbf ' \xcn \xbq ',$ so $ \xbf \xco \xbf
' \xcn \xbq '.$ The following example shows
that this is not always justified.

\er

\be

$\hspace{0.01em}$


\label{Example Drugs}

Consider the following obligations for a physician:

Let $ \xbf ' $ imply that the patient has no heart disease, and if $ \xbf
' $ holds,
we should give drug $ \xCf A$ or (not drug $ \xCf A,$ but drug $ \xCf B),$
abbreviated $A \xco ( \xCN A \xcu B).$
$( \xCf B$ is considered dangerous for people with heart problems.)

Let $ \xbf $ imply that the patient has heart problems. Here, the
obligation is
$(A \xco ( \xCN A \xcu B)) \xcu (A \xco ( \xCN A \xcu \xCN B)),$
equivalent to $ \xCf A.$

The false conclusion would then be $ \xbf ' \xcn A \xco ( \xCN A \xcu B),$
and $ \xbf \xcn A \xco ( \xCN A \xcu B),$
so $ \xbf \xco \xbf ' \xcn A \xco ( \xCN A \xcu B),$ so in both situation
we should either give $ \xCf A$ or
$B,$ but $B$ is dangerous in ``one half'' of the situations.

$ \xcz $
\\[3ex]

\ee

We captured this idea about ``and'' and ``or'' in
Definition \ref{Definition ui} (page \pageref{Definition ui}).
\subsection{
Ceteris paribus - a local poperty
}
\label{Section Ceteris}

Basically, the set of points ``in'' an obligation has to be better than the
set
of ``exterior'' points. As
above Example \ref{Example 3-Obligations} (page \pageref{Example 3-Obligations})
with three obligations shows, demanding
that any element inside is better than any element outside, is too strong.
We use instead the ``ceteris paribus'' idea.

``All other things being equal'' seems to play a crucial role in
understanding obligations. Before we try to analyse it, we look for other
concepts which have something to do with it.

The Stalnaker/Lewis semantics for counterfactual conditionals also works
with
some kind of ``ceteris paribus''.
``If it were to rain, $I$ would use an umbrella'' means something like:
``If it were to rain, and there were not a very strong wind'' (there is no
such
wind now), ``if $I$ had an umbrella'' (I have one now), etc., i.e. if things
were mostly as they are now, with the exception that now it does not rain,
and in the situation $I$ speak about it rains, then $I$ will use an
umbrella.

But also theory revision in the AGM sense contains - at least as objective
-
this idea: Change things as little as possible to incorporate some new
information in a consistent way.

When looking at the ``ceteris paribus'' in obligations, a natural
interpretation
is to read it as ``all other obligations being unchanged'' (i.e. satisfied
or not as before). This is then just a Hamming distance considering the
obligations (but not other information).

Then, in particular, if $ \xdo $ is a family of obligations, and
if $x$ and $x' $ are in the same subset $ \xdo ' \xcc \xdo $ of
obligations, then an
obligation derived from $ \xdo $ should not separate them. More precisely,
if
$x \xbe O \xbe \xdo \xcj x' \xbe O \xbe \xdo,$ and $D$ is a derived
obligation, then $x \xbe D \xcj x' \xbe D.$

\be

$\hspace{0.01em}$


\label{Example Spaghetti}

If the only obligation is not to kill, then it should not be
derivable not to kill and to eat spaghetti.

\ee

Often, this is impossible, as obligations are not independent. In this
case,
but also in other situations, we can push ``ceteris paribus'' into an
abstract distance $d$ (as in the Stalnaker/Lewis semantics), which we
postulate
as given, and say that satisfying an obligation makes things better when
going from ``outside'' the obligation to the $d-$closest situation ``inside''.
Conversely, whatever the analysis of ``ceteris paribus'', and given a
quality
order on the situations, we can now define an obligation as a formula
which
(perhaps among other criteria)
``ceteris paribus'' improves the situation when we go from ``outside'' the
formula ``inside''.

A simpler way to capture ``ceteris paribus'' is to connect it directly to
obligations,
see Definition \ref{Definition Delta-O} (page \pageref{Definition Delta-O}).
This is probably too much tied
to independence (see below), and thus too rigid.
\subsection{
Hamming neighbourhoods
}

A combination concept is a Hamming neighbourhood:

$X$ is called a Hamming neighbourhood of the best cases iff for any $x
\xbe X$
and $y$ a best case with minimal distance from $x,$ all elements between
$x$ and $y$ are in $X.$

For this, we need a notion of distance (also to define ``between'' ).
This was made precise
in Definition \ref{Definition Distance} (page \pageref{Definition Distance})
and
Definition \ref{Definition Neighbourhood} (page \pageref{Definition
Neighbourhood}).
\subsection{
Global and mixed global/local properties of obligations
}

We look now at some global properties (or mixtures of global and local)
which
seem desirable for obligations:

 \xEh

 \xDH Downward closure

Consider the following example:

\be

$\hspace{0.01em}$


\label{Example Not-Global}

Let $U':=\{x,x',y,y' \}$ with $x':=pqrs,$ $y':=pqr \xCN s,$ $x:= \xCN
p \xCN qr \xCN s,$ $y:= \xCN p \xCN q \xCN r \xCN s.$

Consider $X:=\{x,x' \}.$

The counting version:

Then $x' $ has quality 4 (the best), $y' $ has quality 3, $x$ has 1, $y$
has 0.

$d_{c}(x',y' )=1,$ $d_{c}(x,y)=1,$ $d_{c}(x,y' )=2.$

\ee

Then above ``ceteris paribus'' criterion is satisfied, as $y' $ and $x$ do
not
``see'' each other, so
$X \xeb_{l,c} \xdC X.$

But $X$ is not downward closed, below $x \xbe X$ is a better element $y'
\xce X.$

This seems an argument against $X$ being an obligation.

The set version:

We still have $x' \xeb_{s}y' \xeb_{s}x \xeb_{s}y.$ As shown in
Fact \ref{Fact Quality-Distance} (page \pageref{Fact Quality-Distance}),
$d_{s}(x,y)$ (and also
$d_{s}(x',y' ))$ and $d_{s}(x,y' )$ are not comparable, so our argument
collapses.

As a matter of fact, we have the result that the ``ceteris paribus''
criterion entails downward closure in the set variant, see
Fact \ref{Fact Local-Global} (page \pageref{Fact Local-Global}).

$ \xcz $
\\[3ex]

Note that a sufficiently rich domain (put elements between $y' $ and $x)$
will make
this local condition (for $ \xeb )$ a global one, so we have here a domain
problem. Domain problems are discussed e.g. in
 \cite{Sch04} and  \cite{GS08a}.

 \xDH Best states

It seems also reasonable to postulate that obligations contain all best
states. In particular, obligations have then to be consistent - under the
condition that best states exist. We are aware that this point can be
debated,
there is, of course, an easy technical way out: we take, when necessary,
unions of obligations to cover the set of ideal cases.
So obligations will be certain ``neighbourhoods'' of the ``best'' situations.

We think, that some such notion of neighbourhood is a good
candidate for a semantics:

 \xEI

 \xDH
A system of neighbourhoods is not necessarily closed under supersets.

 \xDH
Obligations express something like an approximation to the ideal case
where all obligations (if possible, or, as many as possible) are
satisfied,
so we try to be close to the ideal. If we satisfy an obligation, we are
(relatively) close, and stay so as long as the obligation is satisfied.

 \xDH
The notion of neighbourhood expresses the idea of being close, and
containing
everything which is sufficiently close.
Behind ``containing everything which is sufficiently close'' is the idea
of being in some sense convex. Thus, ``convex'' or ``between''
is another basic notion
to be investigated. See here also the discussion of ``between''
in  \cite{Sch04}.

 \xEJ

 \xEj
\subsection{
Soft obligations
}

``Soft'' obligations are obligations which have exceptions. Normally,
one is obliged to do $O,$ but there are cases where one is not obliged.
This is like soft rules, as ``Birds fly'' (but penguins do not), where
exceptions are not explicitly mentioned.

The semantic notions of size are very useful here, too. We will content
ourselves that soft obligations satisfy the postulates of usual
obligations
everywhere except on a small set of cases. For instance, a soft obligation
$O$
should be downward closed ``almost'' everywhere, i.e. for a small subset
of pairs $ \xBc a,b \xBe $ in $U \xCK U$ we accept that $a \xeb b,$ $b \xbe O,$
$a
\xce O.$ We transplanted
a suitable and cautious notion of size from the components to the
product in Definition \ref{Definition Size-Product} (page \pageref{Definition
Size-Product}).

When we look at the requirement to contain the best cases, we might have
to
soften this, too. We will admit that a small set of the ideal cases might
be
excluded. Small can be relative to all cases, or only to all ideal cases.

Soft obligations generate an ordering which takes care of exceptions, like
the normality ordering of birds will take care of penguins: within the
set of pengins, non-flying animals are the normal ones.
Based on this ordering, we define ``derived soft obligations'', they may
have (a small set of) exceptions with respect to this ordering.
\subsection{
Overview of different types of obligations
}

 \xEh

 \xDH Hard obligations. They hold without exceptions, as in the Ten
Commandments. You should not kill.

 \xEh

 \xDH
In the simplest case, they apply everywhere and can be combined
arbitrarily,
i.e. for any $ \xdo ' \xcc \xdo $ there is a model where all $O \xbe \xdo
' $ hold, and no $O' \xbe \xdo - \xdo '.$

 \xDH
In a more complicated case, not all combinations are possible. This is the
same as considering just an arbitrary subset of $U$ with the same set $
\xdo $ of
obligations. This case is very similar to the case of conditional
obligations
(which might not be defined outside a subset of $U),$ and we treat them
together.

A good example is the Considerate Assassin:

\be

$\hspace{0.01em}$


\label{Example Considerate-Assassin}

Normally, one should not offer a cigarette to someone, out of respect for
his health. But the considerate assassin might do so nonetheless, on the
cynical reasoning that the victim's health is going to suffer anyway:

(1) One should not kill, $ \xCN k.$

(2) One should not offer cigarettes, $ \xCN o.$

(3) The assassin should offer his victim a cigarette before killing him,
if $k,$ then $o.$

Here, globally, $ \xCN k$ and $ \xCN o$ is best, but among $k-$worlds, $o$
is better than $ \xCN o.$
The model ranking is $ \xCN k \xcu \xCN o \xeb \xCN k \xcu o \xeb k \xcu o
\xeb k \xcu \xCN o.$

\ee

 \xEj

Recall that an obligation for the whole set need not be an obligation for
a
subset any more, as it need not contain all best states. In this case, we
may have to take a union with other obligations.

 \xDH Soft obligations.

Many obligations have exceptions. Consider the following example:

\be

$\hspace{0.01em}$


\label{Example Library}

You are in
a library. Of course, you should not pour water on a book. But if the book
has caught fire, you should pour water on it to prevent worse damage. In
stenographic style these obligations read: ``Do not pour water on books''.
``If a book is on fire, do pour water on it.'' It is like ``birds fly'', but
``penguins do not fly'', ``soft'' or nonmonotonic obligations, which have
exceptions, which are not formulated in the original obligation, but added
as exceptions.

\ee

We could have formulated the library obligation also without exceptions:
``When you are in a library, and the book is not on fire, do not pour water
on
it.''
``When you are in a library, and the book is on fire, pour water on it.''
This formulation avoids exceptions. Conditional obligations behave like
restricted quantifiers: they apply in a subset of all possible cases.

We treat now the considerate assassin case as an obligation
(not to offer) with exceptions.
Consider the full set $U,$ and
consider the obligation $ \xCN o.$ This is not downward
closed, as $k \xcu o$ is better than $k \xcu \xCN o.$
Downward closure will only hold for ``most'' cases, but not for all.

 \xDH

Contrary-to-duty obligations.

Contrary-to-duty obligations are about different degrees of fulfillment.
If you should ideally not have any fence, but are not willing or able to
fulfill this obligation (e.g. you have a dog which might stray), then you
should
at least paint it white to make it less conspicuous.
This is also a conditional obligation.
Conditional, as it specifies what has to be done if there is a fence.
The new aspect in contrary-to-duty obligations is the different degree of
fulfillment.

We will not treat contrary-to-duty obligations here, as they do not seem
to
have any import on our basic ideas and solutions.

 \xDH
A still more complicated case is when the language of obligations is not
uniform, i.e. there are subsets $V \xcc U$ where obligations are defined,
which are not defined in $ \xCf U- \xCf V.$

We will not pursue this case here.

 \xEj
\subsection{
Summary of the philosophical remarks
}

 \xEh

 \xDH It seems justifiable to say that an obligation is satisfied or holds
in a certain situation.

 \xDH Obligations are fundamentally asymmetrical, thus negation has to be
treated with care. ``Or'' and ``and'' behave as for facts.

 \xDH Satisfying obligations improves the situation with respect to some
given
grading - ceteris paribus.

 \xDH ``Ceteris paribus'' can be defined by minimal change with respect to
other
obligations, or by an abstract distance.

 \xDH Conversely, given a grading and some distance, we can define an
obligation
locally as describing an improvement with respect to this grading when
going
from ``outside'' to the closest point ``inside'' the obligation.

 \xDH Obligations should also have global properties: they should be
downward (i.e. under increasing quality) closed, and cover the set
of ideal cases.

 \xDH The properties of ``soft'' obligations, i.e. with exceptions, have to
be
modified appropriately. Soft obligations generate an ordering, which in
turn may generate other obligations, where exceptions to the ordering
are permitted.

 \xDH Quality and distance can be defined from an existing set of
obligations
in the set or the counting variant. Their behaviour is quite different.

 \xDH We distinguished various cases of obligations, soft and hard, with
and without all possibilities, etc.

 \xEj

Finally, we should emphasize that the notions of distance, quality, and
size
are in principle independent, even if they may be based on a common
substructure.
\section{
Examination of the various cases
}
\label{Section Exam}

We will concentrate here on the set version of hard obligations.
\subsection{
Hard obligations for the set approach
}
\label{Section Hard-Obligations}
\subsubsection{
Introduction
}

We work here in the set version, the $ \xbe -$case, and examine mostly
the set version only.

We will assume a set $ \xdo $ of obligations to be given.
We define the relation $ \xeb:= \xeb_{ \xdo }$ as described in
Definition \ref{Definition Distance} (page \pageref{Definition Distance}), and
the distance $d$ is the Hamming
distance
based on $ \xdo.$
\subsubsection{
The not necessarily independent case
}

\be

$\hspace{0.01em}$


\label{Example Dependent-1}

Work in the set variant. We show that $X$ $ \xec_{s}-closed$ does not
necessarily imply that $X$ contains all $ \xec_{s}-best$ elements.

Let $ \xdo:=\{p,q\},$ $U':=\{p \xCN q, \xCN pq\},$ then all elements of
$U' $ have best quality
in $U',$ $X:=\{p \xCN q\}$ is closed, but does not contain all best
elements. $ \xcz $
\\[3ex]

\ee

\be

$\hspace{0.01em}$


\label{Example Dependent-3}

Work in the set variant. We show that $X$ $ \xec_{s}-closed$ does not
necessarily imply that $X$ is a neighbourhood of the best elements, even
if $X$ contains them.

Consider $x:=pq \xCN rstu,$ $x':= \xCN pqrs \xCN t \xCN u,$ $x'':=p \xCN
qr \xCN s \xCN t \xCN u,$
$y:=p \xCN q \xCN r \xCN s \xCN t \xCN u,$ $z:=pq \xCN r \xCN s \xCN t
\xCN u.$ $U:=\{x,x',x'',y,z\},$ the $ \xeb_{s}-$best elements
are $x,x',x'',$ they are contained in $X:=\{x,x',x'',z\}.$
$d_{s}(z,x)=\{s,t,u\},$
$d_{s}(z,x' )=\{p,r,s\},$ $d_{s}(z,x'' )=\{q,r\},$ so $x'' $ is one of the
best
elements closest to $z.$ $d(z,y)=\{q\},$ $d(y,x'' )=\{r\},$ so $[z,x''
]=\{z,y,x'' \},$ $y \xce X,$ but $X$
is downward closed. $ \xcz $
\\[3ex]

\ee

\bfa

$\hspace{0.01em}$


\label{Fact Global-Dependent-S}

Work in the set variant.

Let $X \xEd \xCQ,$ $X$ $ \xec_{s}-closed.$ Then

(1) $X$ does not necessarily contain all best elements.

Assume now that $X$ contains, in addition, all best elements. Then

(2) $X \xeb_{l,s} \xdC X$ does not necessarily hold.

(3) $X$ is (ui).

(4) $X \xbe \xdd ( \xdo )$ does not necessarily hold.

(5) $X$ is not necessarily a neighbourhood of the best elements.

(6) $X$ is an improving neighbourhood of the best elements.

\efa

\subparagraph{
Proof
}

$\hspace{0.01em}$


(1) See Example \ref{Example Dependent-1} (page \pageref{Example Dependent-1})

(2) See Example \ref{Example Dependent-2} (page \pageref{Example Dependent-2})

(3)
If there is $m \xbe X,$ $m \xce O$ for all $O \xbe \xdo,$ then by closure
$X=U,$ take $ \xdo_{i}:= \xCQ.$

For $m \xbe X$ let $ \xdo_{m}:=\{O \xbe \xdo:m \xbe O\}.$ Let $X':= \xcV
\{ \xcS \xdo_{m}:m \xbe X\}.$

$X \xcc X':$ trivial, as $m \xbe X \xcp m \xbe \xcS \xdo_{m} \xcc X'.$

$X' \xcc X:$ Let $m' \xbe \xcS \xdo_{m}$ for some $m \xbe X.$ It suffices
to show that $m' \xec_{s}m.$
$m' \xbe \xcS \xdo_{m}= \xcS \{O \xbe \xdo:m \xbe O\},$ so for all $O
\xbe \xdo $ $(m \xbe O \xcp m' \xbe O).$

(4) Consider Example \ref{Example Dependent-2} (page \pageref{Example
Dependent-2}),
let $dom( \xbd )=\{r,s\},$ $ \xbd (r)= \xbd (s)=0.$ Then $x,y \xcm \xbd,$
but $x' \xcM \xbd $ and $x \xbe X,$ $y \xce X,$
but there is no $z \xbe X,$ $z \xcm \xbd $ and $z \xeb y,$ so $X \xce \xdd
( \xdo ).$

(5) See Example \ref{Example Dependent-3} (page \pageref{Example Dependent-3}).

(6) By Fact \ref{Fact Quality-Distance} (page \pageref{Fact Quality-Distance}),
(5).

$ \xcz $
\\[3ex]

\bfa

$\hspace{0.01em}$


\label{Fact Global-Dependent-Converse-S}

Work in the set variant

(1.1) $X \xeb_{l,s} \xdC X$ implies that $X$ is $ \xec_{s}-closed.$

(1.2) $X \xeb_{l,s} \xdC X$ $ \xch $ $X$ contains all best elements

(2.1) $X$ is (ui) $ \xch $ $X$ is $ \xec_{s}-closed.$

(2.2) $X$ is (ui) does not necessarily imply that $X$ contains all $
\xec_{s}-best$
elements.

(3.1) $X \xbe \xdd ( \xdo )$ $ \xch $ $X$ is $ \xec_{s}-closed$

(3.2) $X \xbe \xdd ( \xdo )$ implies that $X$ contains all $
\xec_{s}-best$ elements.

(4.1) $X$ is an improving neighbourhood of the $ \xec_{s}-best$ elements $
\xch $ $X$ is
$ \xec_{s}-closed.$

(4.2) $X$ is an improving neighbourhood of the best elements $ \xch $ $X$
contains all best elements.

\efa

\subparagraph{
Proof
}

$\hspace{0.01em}$


(1.1) By Fact \ref{Fact Local-Global} (page \pageref{Fact Local-Global}).

(1.2) By Fact \ref{Fact General-Obligation} (page \pageref{Fact
General-Obligation}).

(2.1)
Let $O \xbe \xdo,$ then $O$ is downward closed (no $y \xce O$ can be
better than $x \xbe O).$
The rest follows from Fact \ref{Fact Subset-Closure} (page \pageref{Fact
Subset-Closure})  (3).

(2.2) Consider Example \ref{Example Dependent-1} (page \pageref{Example
Dependent-1}), $p$ is (ui) (formed in
$U!),$ but $p \xcs X$ does
not contain $ \xCN pq.$

(3.1)
Let $X \xbe \xdd ( \xdo ),$ but let $X$ not be closed.
Thus, there are $m \xbe X,$ $m' \xec_{s}m,$ $m' \xce X.$

Case 1: Suppose $m' \xCq m.$ Let $ \xbd_{m}: \xdo \xcp 2,$ $
\xbd_{m}(O)=1$ iff $m \xbe O.$ Then $m,m' \xcm \xbd_{m},$
and there cannot be any $m'' \xcm \xbd_{m},$ $m'' \xeb_{s}m',$ so $X \xce
\xdd ( \xdo ).$

Case 2: $m' \xeb_{s}m.$ Let $ \xdo ':=\{O \xbe \xdo:m \xbe O \xcj m'
\xbe O\},$ $dom( \xbd )= \xdo ',$ $ \xbd (O):=1$ iff $m \xbe O$ for
$O \xbe \xdo '.$ Then $m,m' \xcm \xbd.$ If there is $O \xbe \xdo $ s.t.
$m' \xce O,$ then by $m' \xec_{s}m$ $m \xce O,$
so $O \xbe \xdo '.$ Thus for all $O \xce dom( \xbd ).m' \xbe O.$ But then
there is no $m'' \xcm \xbd,$ $m'' \xeb_{s}m',$
as $m' $ is already optimal among the $n$ with $n \xcm \xbd.$

(3.2) Suppose $X \xbe \xdd ( \xdo ),$ $x' \xbe U-X$ is a best element,
take $ \xbd:= \xCQ,$ $x \xbe X.$
Then there must be $x'' \xeb x',$ $x'' \xbe X,$ but this is impossible as
$x' $ was best.

(4.1) By Fact \ref{Fact Quality-Distance} (page \pageref{Fact Quality-Distance})
, (4) all minimal elements
have incomparabel
distance. But if $z \xec y,$ $y \xbe X,$ then either $z$ is
minimal or it is above a minimal element, with minimal distance from $y,$
so $z \xbe X$
by Fact \ref{Fact Quality-Distance} (page \pageref{Fact Quality-Distance})  (3).

(4.2) Trivial.

$ \xcz $
\\[3ex]
\subsubsection{
The independent case
}

Assume now the system to be independent, i.e. all combinations of $ \xdo $
are present.

Note that there is now only one minimal element, and the notions of
Hamming neighbourhood of the best elements
and improving Hamming neighbourhood of the best elements coincide.

\bfa

$\hspace{0.01em}$


\label{Fact Global-Independent-S}

Work in the set variant.

Let $X \xEd \xCQ,$ $X$ $ \xec_{s}-closed.$ Then

(1) $X$ contains the best element.

(2) $X \xeb_{l,s} \xdC X$

(3) $X$ is (ui).

(4) $X \xbe \xdd ( \xdo )$

(5) $X$ is a (improving) Hamming neighbourhood of the best elements.

\efa

\subparagraph{
Proof
}

$\hspace{0.01em}$


(1) Trivial.

(2) Fix $x \xbe X,$ let $y$ be closest to $x,$ $y \xce X.$ Suppose $x \xeB
y,$ then there must be
$O \xbe \xdo $ s.t. $y \xbe O,$ $x \xce O.$ Choose $y' $ s.t. $y' $ is
like $y,$ only $y' \xce O.$ If $y' \xbe X,$ then
by closure $y \xbe X,$ so $y' \xce X.$ But $y' $ is closer to $x$ than $y$
is, $contradiction.$

Fix $y \xbe U-$X. Let $x$ be closest to $y,$ $x \xbe X.$ Suppose $x \xeB
y,$ then there is $O \xbe \xdo $
s.t. $y \xbe O,$ $x \xce O.$ Choose $x' $ s.t. $x' $ is like $x,$ only $x'
\xbe O.$ By closure of $X,$
$x' \xbe X,$ but $x' $ is closer to $y$ than $x$ is, $contradiction.$

(3) By Fact \ref{Fact Global-Dependent-S} (page \pageref{Fact
Global-Dependent-S})  (3)

(4)
Let $X$ be closed, and $ \xdo ' \xcc \xdo,$ $ \xbd: \xdo ' \xcp 2,$
$m,m' \xcm \xbd,$ $m \xbe X,$ $m' \xce X.$
Let $m'' $ be s.t. $m'' \xcm \xbd,$ and for all $O \xbe \xdo -dom( \xbd
)$ $m'' \xbe O.$ This exists by
independence.
Then $m'' \xec_{s}m',$ but also $m'' \xec_{s}m,$ so $m'' \xbe X.$ Suppose
$m'' \xCq m',$ then
$m' \xec_{s}m'',$ so $m' \xbe X,$ contradiction, so $m'' \xeb_{s}m'.$

(5) Trivial by (1), the remark preceding this Fact, and
Fact \ref{Fact Global-Dependent-S} (page \pageref{Fact Global-Dependent-S})
(6).

\bfa

$\hspace{0.01em}$


\label{Fact Global-Independent-Converse-S}

Work in the set variant.

(1) $X \xeb_{l,s} \xdC x$ $ \xch $ $X$ is $ \xec_{s}-closed,$

(2) $X$ is (ui) $ \xch $ $X$ is $ \xec_{s}-closed,$

(3) $X \xbe \xdd ( \xdo )$ $ \xch $ $X$ is $ \xec_{s}-closed,$

(4) $X$ is a (improving) neighbourhood of the best elements $ \xch $ $X$
is $ \xec_{s}-closed.$

\efa

\subparagraph{
Proof
}

$\hspace{0.01em}$


(1) Suppose there are $x \xbe X,$ $y \xbe U-$X, $y \xeb x.$ Choose them
with minimal distance.
If $card(d_{s}(x,y))>1,$ then there is
$z,$ $y \xeb_{s}z \xeb_{s}x,$ $z \xbe X$ or $z \xbe U-$X, contradicting
minimality. So $card(d_{s}(x,y))=1.$
So $y$ is among the closest elements of $U-X$ seen from $x,$ but then by
prerequisite
$x \xeb y,$ $contradiction.$

(2) By Fact \ref{Fact Global-Dependent-Converse-S} (page \pageref{Fact
Global-Dependent-Converse-S})  (2.1).

(3) By Fact \ref{Fact Global-Dependent-Converse-S} (page \pageref{Fact
Global-Dependent-Converse-S})  (3.1).

(4) There is just one best element $z,$ so if $x \xbe X,$ then $[x,z]$
contains all
$y$ $y \xeb x$ by Fact \ref{Fact Quality-Distance} (page \pageref{Fact
Quality-Distance})  (3).

$ \xcz $
\\[3ex]

The $ \xdd ( \xdo )$ condition seems to be adequate only for the
independent situation,
so we stop considering it now.

\bfa

$\hspace{0.01em}$


\label{Fact Int-Union}

Let $X_{i} \xcc U,$ $i \xbe I$ a family of sets, we note the following
about closure under unions and intersections:

(1) If the $X_{i}$ are downward closed, then so are their unions and
intersections.

(2) If the $X_{i}$ are (ui), then so are their unions and intersections.

\efa

\subparagraph{
Proof
}

$\hspace{0.01em}$


Trivial. $ \xcz $
\\[3ex]

We do not know whether $ \xeb_{l,s}$ is preserved under unions and
intersections,
it does not seem an easy problem.

\bfa

$\hspace{0.01em}$


\label{Fact Relativization}

(1) Being downward closed is preserved while going to subsets.

(2) Containing the best elements is not preserved (and thus neither the
neighbourhood property).

(3) The $ \xdd ( \xdo )$ property is not preserved.

(4) $ \xec_{l,s}$ is not preserved.

\efa

\subparagraph{
Proof
}

$\hspace{0.01em}$


(4) Consider Example \ref{Example Not-Global} (page \pageref{Example
Not-Global}), and
eliminate $y$ from $U',$ then the closest to $x$ not in $X$ is $y',$
which is better.

$ \xcz $
\\[3ex]
\subsection{
Remarks on the counting case
}

\br

$\hspace{0.01em}$


\label{Remark}

In the counting variant all qualities are comparabel. So if $X$
is closed, it will contain all minimal elements.

\er

\be

$\hspace{0.01em}$


\label{Example H-N-Local}

We measure distance by counting.

Consider $a:= \xCN p \xCN q \xCN r \xCN s,$ $b:= \xCN p \xCN q \xCN rs,$
$c:= \xCN p \xCN qr \xCN s,$ $d:=pqr \xCN s,$
let $U:=\{a,b,c,d\},$ $X:=\{a,c,d\}.$ $d$ is the best element,
$[a,d]=\{a,d,c\},$
so $X$ is an improving Hamming neighbourhood, but $b \xeb a,$ so $X
\xeB_{l,c} \xdC X.$

$ \xcz $
\\[3ex]

\ee

\bfa

$\hspace{0.01em}$


\label{Fact Local-H-N}

We measure distances by counting.

$X \xeb_{l,c} \xdC X$ does not necessarily imply that $X$ is an improving
Hamming
neighbourhood of the best elements.

\efa

\subparagraph{
Proof
}

$\hspace{0.01em}$


Consider Example \ref{Example Not-Global} (page \pageref{Example Not-Global}).
There $X \xeb_{l,c} \xdC X.$
$x' $ is the best element, and
$y' \xbe [x',x],$ but $y' \xce X.$ $ \xcz $
\\[3ex]
\section{
What is an obligation?
}
\label{Section What}

The reader will probably not expect a final definition. All we can do is
to give a tentative definition, which, in all probability, will not be
satisfactory in all cases.

\bd

$\hspace{0.01em}$


\label{Definition Obligation}

We decide for the set relation and distance.

(1) Hard obligation

A hard obligation has the following properties:

(1.1) It contains all ideal cases in the set considered.

\ed

(1.2) It is closed under increasing quality,
Definition \ref{Definition Closed} (page \pageref{Definition Closed})

(1.3) It is an improving neighbourhood of the ideal cases (this also
implies (1.1)),
Definition \ref{Definition Neighbourhood} (page \pageref{Definition
Neighbourhood})

We are less committed to:

(1.4) It is ceteris paribus improving,
Definition \ref{Definition Quality-Extension} (page \pageref{Definition
Quality-Extension})

An obligation $O$ is a derived obligation of a system $ \xdo $ of
obligations iff
it is a hard obligation based on the set variant of the order and distance
generated by $ \xdo.$

(2) Soft obligations

A set is a soft obligation iff it satisfies the soft versions of
above postulates. The notion of size has to be given, and is transferred
to products as described in Definition \ref{Definition Size-Product} (page
\pageref{Definition Size-Product}).
More precisely, strict universal quantifiers are transformed into their
soft variant ``almost all'', and the other operators are left as they are.
Of course, one might also want to use a mixture of soft and hard
conditions, e.g. we might want to have all ideal cases, but renounce
on closure for a small set of pairs $ \xBc x,x'  \xBe.$

An obligation $O$ is derived from $ \xdo $ iff it is a soft obligation
based
on the set variant of the order and distance generated by the
translation of $ \xdo $ into their hard versions. (I.e. exceptions will
be made explicit.)

\bfa

$\hspace{0.01em}$


\label{Fact Garbage-In}

Let $O \xbe \xdo,$ then $ \xdo \xcn O$ in the independent set case.

\efa

\subparagraph{
Proof
}

$\hspace{0.01em}$


We check $(1.1)-(1.3)$ of Definition \ref{Definition Obligation} (page
\pageref{Definition Obligation}).

(1.1) holds by independence.

(1.2) If $x \xbe O,$ $x' \xce O,$ then $x' \xeC_{s}x.$

(1.3) By Fact \ref{Fact Global-Dependent-S} (page \pageref{Fact
Global-Dependent-S})  (6).

Note that (1.4) will also hold by
Fact \ref{Fact Global-Independent-S} (page \pageref{Fact Global-Independent-S})
(2).

$ \xcz $
\\[3ex]

\bco

$\hspace{0.01em}$


\label{Corollary Classical-Consequence}

Every derived obligation is a classical consequence of the original set
of obligations in the independent set case.

\eco

\subparagraph{
Proof
}

$\hspace{0.01em}$


This follows
from Fact \ref{Fact Global-Independent-S} (page \pageref{Fact
Global-Independent-S})  (3)
and Fact \ref{Fact Garbage-In} (page \pageref{Fact Garbage-In}).

\be

$\hspace{0.01em}$


\label{Example Burnt-Letter}

The Ross paradox is not a derived obligation.

\ee

\subparagraph{
Proof
}

$\hspace{0.01em}$


Suppose we have the alphabet $p,q$ and the obligations $\{p,q\},$ let
$R:=p \xco \xCN q.$ This is not not closed, as $ \xCN p \xcu q \xeb \xCN p
\xcu \xCN q \xbe R.$
$ \xcz $
\\[3ex]
\section{
Conclusion
}

Obligations differ from facts in the behaviour of negation, but not of
conjunction and disjunction. The Ross paradox originates, in our opinion,
from the differences in negation. Central to the treatment of obligations
seems to be a relation of ``better'', which can generate obligations, but
also be generated by obligations. The connection between obligations and
this relation of ``better'' seems to be somewhat complicated and leads to
a number of ramifications. A tentative definition of a derivation of
obligations is given.
\section{
Acknowledgements
}

We thank A.Herzig, Toulouse, and L.v.d.Torre, Luxembourg, for very helpful
discussions.

\end{document}